\documentclass{article}   

\usepackage{amsmath,amssymb}
\usepackage{amsthm}

\title{b-coloring graphs with large girth}
\author{%
Victor Campos \thanks{ParGO, Universidade Federal do Cear\'a, Brazil} \thanks{E-mail: campos@lia.ufc.br, victorfarias@lia.ufc.br, anasilva@mat.ufc.br} %
\and %
Victor Farias \footnotemark[1] \footnotemark[2] %
\and %
Ana Silva \footnotemark[1] \footnotemark[2]%
}

\newtheorem{theorem}{Theorem}[section]

\begin{document}

%\begin{frontmatter}

% \title{b-coloring graphs with large girth\tnoteref{ALL}}
% \tnotetext[ALL]{Partially supported by Funcap and CNPq/ Brazil.}
% 
% \author[Comp]{Victor Aguiar Evangelista de Farias}
% \ead{victorfarias@lia.ufc.br}
% 
% \author[UFCVirtual]{Victor Campos}
% \ead{campos@lia.ufc.br}
% 
% \author[Mat]{Ana Silva}
% \ead{ana.silva@g-scop.inpg.fr}
% 
% \address[Comp]{Departamento de Computa\c{c}\~ao, Universidade Federal do Cear\'a}
% \address[UFCVirtual]{Instituto UFC Virtual, Universidade Federal do Cear\'a}
% \address[Mat]{Departamento de Matem\'atica, Universidade Federal do Cear\'a}

%\author{Victor Campos \and Victor Farias\and Ana Silva}
%
%\institute{\at Universidade Federal do Cear\'a, Brazil. \email{$\{$campos,victorfarias$\}$@lia.ufc.br} and anasilva@mat.ufc.br.}
% \ead{victorfarias@lia.ufc.br}
% 
% \author[UFCVirtual]{Victor Campos}
% \ead{campos@lia.ufc.br}
% 
% \author[Mat]{Ana Silva}
% \ead{ana.silva@g-scop.inpg.fr}
% 
% \address[Comp]{Departamento de Computa\c{c}\~ao, Universidade Federal do Cear\'a}
% \address[UFCVirtual]{Instituto UFC Virtual, Universidade Federal do Cear\'a}
% \address[Mat]{Departamento de Matem\'atica, Universidade Federal do Cear\'a}


\maketitle

\begin{abstract}
A b-coloring of a graph is a coloring of its vertices such that every
color class contains a vertex that has a neighbor in all other
classes.  The b-chromatic number of a graph is the largest integer $k$
such that the graph has a b-coloring with $k$ colors.  We show
how to compute in polynomial time the b-chromatic number of a graph of girth at least $9$.  This improves the seminal
result of Irving and Manlove on trees.
%\keywords{b-chromatic number \and b-coloring \and m-degree \and girth \and exact algorithm}
\end{abstract}

\section{Introduction}

Let $G$ be a simple graph.  A \textit{proper coloring of $G$} is an
assignment of colors to the vertices of $G$ such that no two adjacent
vertices have the same color.  The \textit{chromatic number of $G$}
is the minimum integer $\chi(G)$ such that $G$ has a proper coloring
with $\chi(G)$ colors.  Suppose that we have a proper coloring of
$G$ and there exists a color $h$ such that every vertex $x$ with
color $h$ is not adjacent to at least one other color (which may
depend on $x$); then we can change the color of these vertices and
thus obtain a proper coloring with fewer colors.  This heuristic can
be applied iteratively, but we cannot expect to reach the chromatic
number of $G$, since the coloring problem is $\mathcal{NP}$-hard.
On the basis of this idea, Irving and Manlove introduced the notion of
b-coloring in \cite{Irving.Manlove.99}.  
Intuitively, a b-coloring is a proper coloring that cannot be improved by the above heuristic, and the b-chromatic number measures the worst possible such coloring. 
More formally, consider any vertex coloring of $G$.  
A vertex $u$ is said to be a \emph{b-vertex} (for this coloring) if $u$ has a neighbor colored with each color different from its own color.  
A \emph{b-coloring} of $G$ is a proper coloring of $G$ such that each color class contains a b-vertex.  
A {\em basis} of a b-coloring is a set of b-vertices, one for each color class.
The \emph{b-chromatic number} of $G$ is the largest integer $k$ such that $G$ has a b-coloring with $k$ colors.  
We denote it by $\chi_b(G)$.
%In a b-coloring with $k$ colors, let $v_i$ be any b-vertex of color $i$ ($i=1, \ldots, k$); then we say that the set $\{v_1, \ldots,v_k\}$ is a \emph{basis} of the b-coloring.  A b-coloring may have many bases.

Naturally, a proper coloring of $G$ with $\chi(G)$ colors is a b-coloring of $G$, since it cannot be improved.  
Hence, $\chi(G)\leq \chi_b(G)$.  
For an upper bound, observe that if $G$ has a b-coloring with $k$ colors, then $G$ has at least $k$ vertices with degree at least $k-1$ (a basis of the b-coloring). 
Thus, if $m(G)$ is the largest integer such that $G$ has at least $m(G)$ vertices with degree at least $m(G)-1$, we know that $G$ cannot have a b-coloring with more than $m(G)$ colors, i.e., % 
$$\chi_b(G)\leq m(G).$$
This upper bound was introduced by Irving and Manlove in \cite{Irving.Manlove.99}.  
They showed that the difference between $\chi_b(G)$ and $m(G)$ can be arbitrarily large for general graphs.
They proved that $\chi_b(G)$ is equal to $m(G)$ or $m(G)-1$ when $G$ is a tree, and provided a polynomial time algorithm that computes $\chi_b(G)$ for every tree.  
In addition, the problem was proved to be NP-hard in general graphs \cite{Irving.Manlove.99}, and remains so even when restricted to bipartite graphs \cite{KRATOCHVIL.etal.02}. 
These concepts have received much attention recently; for example, see \cite{Balakrishnan.Raj.09} to \cite{Velasquez.Bonomo.Koch.11}.

Many of these works investigate the b-chromatic number of graphs under assumptions that involve the existence of large cycles. For example, Irving and Manlove's algorithm for trees can actually work on graphs with girth at least 11, as noticed by A.~Silva in \cite{Silva.10}. Also, there are a number of results about d-regular graphs with girth at least~5 \cite{Blidia.etal.09,Cabello.Jakovac.10,Kouider.04,KRATOCHVIL.etal.02,Sahili.Kouider.06}. In this paper we improve Irving and Manlove's result for graphs with large girth; more specifically, we prove the following.

\begin{theorem}
 If $G$ is a graph with girth at least 9, then $\chi_b(G)\geq m(G)-1$.
\label{thm:main}
\end{theorem}

Here is an outline of the proof of Theorem~\ref{thm:main}.  A special set of vertices, called a good set of vertices, is defined and graphs are distinguished between having a good set and not having a good set. Next, we state some results by Irving and Manlove \cite{Irving.Manlove.99} and by A.~Silva \cite{Silva.10} that say that a graph $G$ with $\mbox{girth}(G)\geq 8$ that does not have a good set cannot be b-colored with $m(G)$ colors and has a b-coloring with $m(G)-1$ colors (hence, $\chi_b(G)=m(G)-1$); also, A.~Silva proved that if $G$ with girth at least 8 has a good set, then one can be found in polynomial time. Finally, and this is the original part of the paper, it is shown that if $G$ with girth at least 9 has a good set, then $\chi_b(G)=m(G)$.
The proof of Theorem~\ref{thm:main} yields a polynomial time algorithm that finds an optimal b-coloring of graphs with girth at least 9.
 %Actually, a part of this proof has already been done by , as we will see in Section \ref{sec:previousResults}.

% Consider $G$ to be a graph with girth at least 9. In Section \ref{sec:defs}, we give the basic definitions and present the notation that will be used along the text. In Section \ref{sec:previousResults}, we present the definitions of good set and the results by A.~Silva \cite{Silva.10} that proves that a good set of $G$ can be found in polynomial time, if one exists, and that if one does not exist, then $\chi_b(G)=m(G)-1$. Finally, in Section \ref{sec:nonpivotedColouring}, we show that if $G$ has a good set, then $\chi_b(G)=m(G)$. 
% 
% \section{Definitions and notation}
% \label{sec:defs}

%%%%%
\section{Definitions and partial results}
\label{sec:previousResults}

In this section, we present some necessary definitions and the results by Irving and Manlove \cite{Irving.Manlove.99} and A.~Silva \cite{Silva.10} that complement our proof. The graph terminology used in this paper follows \cite{BM08}. 

Let $G$ be a simple graph. We denote by $V(G)$ and $E(G)$ the sets of vertices and edges of $G$, respectively.  If $X\subseteq V(G)$, then $N^X(u)$ represents the set $N(u)\cap X$. The \emph{girth of $G$} is the size of a shortest induced cycle of $G$.

Recall that $m(G)$ is the largest integer $k$ such that $G$ has at
least $k$ vertices with degree at least $k-1$.  We say that a vertex
$u\in V(G)$ is \emph{dense} if $d(u)\geq m(G)-1$; and we denote the set
of dense vertices of $G$ by $M(G)$.

Let $W$ be a subset of $M(G)$, and let $u$ be any vertex in $V(G)\setminus W$.  If $u$ is such that every vertex $v\in W$ is either adjacent to $u$ or has a common neighbor $w\in W$ with $u$ such that $d(w)=m(G)-1$, then it is said that \textit{$W$ encircles vertex $u$} (or that $u$ is encircled by $W$). A subset $W$ of $M(G)$ of size $m(G)$ is a \emph{good set} if (our definition is slightly different from the one given by Irving and Manlove): \\
(a) $W$ does not encircle any vertex, and \\
(b) Every vertex $x\in V(G)\setminus W$ with $d(x) \ge m(G)$ is adjacent to a vertex $w\in W$.% with $d(w) = m(G)-1$.

\begin{theorem}[\cite{Irving.Manlove.99}]
 Let $G$ be any graph and $W$ be a subset of $M(G)$ with $m(G)$ vertices. If $W$ encircles some vertex $v\in V(G)\setminus W$, then $W$ is not a basis of a b-coloring with $m(G)$ colors.
\end{theorem}

\begin{theorem}[\cite{Silva.10}]
If $G$ is a graph with girth at least~$8$, then $G$ does not have a good set if and only if $|M(G)|=m(G)$ and $M(G)$ encircles a vertex in $V(G)\setminus M(G)$.  Moreover, a good set of $G$ (if any exists) can be found in polynomial time.
\label{thm:findGoodset}
\end{theorem}

A part of the proof of Theorem \ref{thm:main} consists of the following theorem:

\begin{theorem}[\cite{Silva.10}]
\label{thm:nogoodset}
 Let $G$ be a graph with girth at least~$8$. If $G$ has no good set, then $\chi_b(G)=m(G)-1$.
\end{theorem}

Now, all we need to prove is that if $G$ does have a good set, then $G$ can be b-colored with $m(G)$ colors, which is done in the next section.

%%%%%
\section{Coloring graphs with a good set}
\label{sec:nonpivotedColouring}

In this section we prove the second part of the main theorem, namely:
\begin{theorem}\label{thm:goodset}
Let $G$ be a graph with girth at least~$9$.  If $G$ has a good set, then $\chi_b(G)=m(G)$.
\end{theorem}

Let $W = \{v_1, \ldots, v_{m(G)}\}$ be a good set of $G$.  Our aim is
to construct a b-coloring of $G$ with $m(G)$ colors such that, for
each $i\in\{1, \ldots, m(G)\}$, vertex $v_i$ is a b-vertex of color
$i$.  We start by assigning color $i$ to $v_i$, for each $i\in\{1,\ldots, m(G)\}$. Next, we extend this partial coloring to the rest of the
graph in several steps.  Before explaining each step, we need to
introduce some other terminology and notation.

A \emph{link} is any path of length two or three whose extremities are in $W$ and whose internal vertices are not in $W$.  Any interior vertex of a link is called a \emph{link vertex}.  Let $L$ be the set of all link vertices. %For every $w\in W$ and vertex $x\in L$ adjacent to $w$, let $D(w,x)$ be the set of colors $j$ such that $x$ lies on a link between $w$ and $v_j$, i.e., on a link $\langle w, x, v_j\rangle$ or $\langle w,x, x', v_j\rangle$.

We first color $G[W\cup L]$ in a way not to repeat too many colors in $N(w)$, for all $w\in W$, and at the end we extend the obtained partial coloring to a b-coloring of $G$ with $m(G)$ colors. Let $G'=G[W\cup L]$, $L_1$ be the set of vertices of $L$ that have at least one neigbour in $L$ and $L_2$ be the set of vertices in $L$ that have at least two neigbours in $W$. The steps below are followed in order in such a way that we only move on to the next step when all the possible vertices are iterated.

\begin{enumerate}
 \item For each $x\in L_1$, let $x'\in N^L(x)$. Since $x'\in L$, there must exist $v_i\in N^W(x')$; color $x$ with $i$;
 \item For each $v_i\in W$, let $N^*_i=N(v_i)\cap L_2=\{x_1,\ldots,x_q\}$. Also, let $v_{i_j}\in N^W(x_j)\setminus \{v_i\}$. If $q>1$, then use colors $i_1,\ldots,i_q$ to color the uncolored vertices in $N^*_i$ in a way that $x_j$ is not colored with $i_j$ (it suffices to make a derangement of those colors on the vertices);
 \item Let $x\in L_2$ still uncolored be such that there exists $v_i\in N^W(x)$ that has some neighbor $y\in L_1$. Let $c$ be the color of $y$; color $x$ with $c$ and recolor $y$ with $j$, for any $v_j\in N^W(x)\setminus\{v_i\}$;
 %Let $x\in L_2$ still uncolored and let $N_x=\{v_i\in N^W(x)\mid N^{L_1}(v_i)\neq\emptyset\}$. Suppose that $N_x\neq \emptyset$ and let $v_i\in N_x$ and $y\in N^{L_1}(v_i)$. Let $c=c(y)$ and let $Y=\{y=y_0,\ldots,y_p\}$ be the set of all vertices in $L_1$ colored with $c$ that have a common neighbor with $x$ in $W$. Note that $\{v_k\mid v_k\in N^W(y_j),\mbox{ for some $y_j\in Y$}\}\subseteq N_x$, since $Y\cap L_2=\emptyset$. Denote by $v_{i_j}$ the vertex in $N^W(x)\cap N^W(y_j)$, for each $j\in\{0,\ldots,p\}$. If $p>0$, color $x$ with $c$ and use the colors $\{i_0,\ldots,i_p\}$ to recolor the vertices in $Y$ in such a way that $y_j$ is not colored with $i_j$. Otherwise, color $x$ with $i$ and recolor $y$ with $j$, for any $v_j\in N^W(x)\setminus\{v_i\}$;
 \item Finally, if $x\in L_2$ is still uncolored, we know that $N^L(v_i)=\{x\}$, for all $v_i\in N^W(x)$. Since $N^L(x)=\emptyset$, we can color $x$ with $i$, for any $v_i$ that is not adjacent to $x$ and has no common neighbor with $x$ in $W$ of degree $m(G)-1$, which exists as $x$ is not encircled by $W$.
\end{enumerate}

Suppose that the algorithm above produces a partial coloring that colors every vertex in $L$ in such a way that, at the end, each  $v_i\in W$ has at least as many uncolored neighbors as missing colors in its neighborhood. Since $L$ is colored, we know that the uncolored neighbors of $W$ form a stable set. Thus, we can independently color $N(v_i)$ in such a way that $v_i$ sees every other color, for all $v_i\in W$. By the definition of a good set, we know that if $d(v)\ge m(G)$, then $v$ is already colored; hence, the partial coloring can be greedily transformed 
into a b-coloring with $m(G)$ colors. Now, to prove that the algorithm works, we show that after these steps the obtained partial coloring $\psi$ satisfies:

\begin{itemize}
 \item[P1] $\psi$ is proper; and 
 \item[P2] the number of uncolored neighbors of $v_i$ is at least the number of missing colors in $N(v_i)$, for each $v_i\in W$.
\end{itemize}

% However, P2 does not always hold after Step 2. Consider, for instance, a cycle $\langle v_{i_1},x_1,v_{i_2},x_2,v_{i_3},x_3,v_{i_4},x_4\rangle$, where $x_i\in L_2$ is not colored, for all $i\in \{1,\ldots,4\}$. If we color first the neigbours of $v_{i_1}$ and then the neighbors of $v_{i_3}$ we might end with a coloring where $x_1$ and $x_2$ are colored with $i_3$ and $x_3$ and $x_4$ are colored with $i_2$. So, we will suppose that we have a partial order $\prec$ of $W$ such that if such a cycle exists and $v_{i_1}\prec v_{i_j}$, for $j=2$ and $j=4$, then $v_{i_j}\prec v_{i_3}$, for $j=2$ or $j=4$. Then, we iterate on $W$ in Step 2 according to this order.
 
\textit{Proof of Theorem \ref{thm:goodset}}: 
% We actually prove that, until Step 3, Property P1 holds and no color is repeated in $N(v_i)$, for each $v_i\in W$, which trivially implies P2. 
% Consider the same notation as above. 

First, we make some observations concerning the coloring procedure. Note that $L_1\cap L_2$ is not necessarily empty, but all vertices in this subset are colored in Step 1. However, a vertex $x\in L_1\cap L_2$ may play a role in Step 2 in the following way: if $x\in N(v_i)$ and there exists $x'\in N^{L_2}(v_i)\setminus L_1$, then $x'$ may be colored with color $j$ for some $v_j\in N^W(x)\setminus \{v_i\}$, while the color of $x$ remains unchanged. Also, note that, in Step 3, since $N^{L_2}(v_i)=\{x\}$, we have $y\in L_1\setminus L_2$. Hence, $N^W(y)=\{v_i\}$ and, consequently, the color of $y$ cannot be changed again. Thus (*) the color of $y$ is changed at most once, for every $y\in L_1$. Finally, if $x$ receives color $i$ in Step~1,~2 or~3, then one of the following holds (fact (iii) holds because of (*)):

\begin{enumerate}
 \item[(i)] {\it $x$ receives color $i$ in Step 1 and there exists a path $\langle x,x',v_i\rangle$, for some $x'\in L_1$; or}
 \item[(ii)] {\it $x$ receives color $i$ in Step 2 and there exists a path $\langle x,v_j,x',v_i\rangle$, for some $v_j\in W$ and $x'\in L_2$; or}
 \item[(iii)] {\it $x$ receives color $i$ in Step 3 and there exists a path $\langle x,v_j,y,y',v_i\rangle$, for some $v_j\in W$, $y\in L_1\setminus L_2$ and $y'\in L_1$; or}
 \item[(iv)] {\it $x$ is recolored with color $i$ in Step 3 and there exists a path $\langle x,v_j,x',v_i\rangle$, for some $v_j\in W$ and $x'\in L_2\setminus L_1$.}
\end{enumerate}

We first prove that P1 holds after Step 3. 
Suppose that there exists an edge $wz$ such that $\psi(w)=\psi(z)=i$. 
Since $G$ has no cycle of length at most 7, the paths defined in (i)-(iv) are shortest paths.
Therefore, vertex $v_i$ has no neighbor colored $i$ and hence, $w,z\in L$. 
Also, as $wz\in E(G)$, we have $w,z\in L_1$ and they are colored in Step 1 and maybe recolored in Step 3. 
By (i) and (iv), there exist a $w,v_i$-path $P_w$ and a $z,v_i$-path $P_{z}$, both of length at most 3. 
Note that either $P_w+P_{z}+wz$ contains a cycle of length at most 7 or one of these paths consists of the edge $wz$ followed by the other path.
Because $G$ has girth at least 9, the latter case occurs.
We get as contradiction as this implies that at least one path is defined by (i) and, thus, vertex $v_i$ has a neighbor colored $i$.
% Also by (i) and (iv), one can see that these paths cannot intersect, a contradiction as in this case we get a cycle of length at most 7.

Now, we prove that P2 also holds after Step 3. 
We actually prove that, after Step 3, no color is repeated in $N(v_i)$, for each $v_i\in W$.
Suppose there exist a vertex $v_j\in W$ and $w,z\in N(v_j)$ such that $\psi(w) = \psi(z) = i$.
First, consider the case $v_i\in \{w,z\}$.
Since the paths defined by (i)-(iv) are shortest paths, we have that (i) occurs for the vertex in $\{w, z\}\setminus \{v_i\}$.
We get a contradiction as this implies $G$ has a cycle of length 4. Therefore we may assume $v_i\notin \{w,z\}$.

Now, by (i)-(iv), there exist a $w,v_i$-path $P_w$ and a $z,v_i$-path $P_{z}$.
Let $\ell_w$ and $\ell_z$ be the length of $P_w$ and $P_z$, respectively. Clearly $\ell_w,\ell_z\le 4$. 
Note that either $P_w+P_{z}+\langle w,v_j,z\rangle$ contains a cycle of length at most $\ell_w + \ell_z + 2$ or either $P_w$ or $P_z$ consists of the path $\langle w, v_j, z\rangle$ followed by the other path.
Since both $w$ and $z$ are at distance at least 2 from $v_i$ and $\ell_w,\ell_z\le 4$, the latter can only occur if one of the paths is defined by (i), say $P_w$, and the other is defined by (iii), say $P_z$.
We get a contradiction as $P_z = \langle z, v_j, w, y, v_i \rangle$ implies $w$ is recolored in Step~3 and therefore, $P_w$ must be defined by (iv). 
Now, suppose that the former occurs, i.e., $P_w+P_{z}+\langle w,v_j,z\rangle$ contains a cycle of length at most $\ell_w + \ell_z + 2$.
Because $G$ has girth at least 9, we have $\ell_w + \ell_z \ge 7$. 
This implies that at least one of $P_w$ and $P_z$, say $P_z$, is defined by (iii), and the other is not defined by (i).
Therefore $z$ is colored in Step 3 and $N^{L_2}(v_j) = \{z\}$.
Furthermore, $w\in L_1\setminus L_2$ and $N^W(w) = \{v_j\}$. Therefore, since (i) does not occur for $w$, we have that $P_w$ must be defined by (iv). Thus the only choice for $P_w$ is $\langle w, v_j, z, v_i\rangle$, a contradiction as P1 holds.

Finally, consider $x$ to be colored during Step 4 with color $i$. By the choice of $i$ we know that $v_i\notin N(x)$. Thus, since $N^L(x)=\emptyset$, Property P1 holds. Now, suppose that some $v_j\in N(x)$ is such that color $i$ already appears in $N(v_j)$. Since $N^L(v_j)=\{x\}$ we must have $v_i\in N(v_j)$ and, by the choice of $i$, $d(v_j)>m(G)-1$. Property P2 thus follows as $i$ is the only repeated color in the neighborhood of $v_j$.\hfill $\Box$

\section{Conclusion}

We showed that if $G$ is a graph with girth at least~$9$, then $\chi_b(G)\geq m(G)-1$, improving the result by Irving and Manlove \cite{Irving.Manlove.99}. We also give an algorithm that finds the b-chromatic number of $G$ in polynomial time. %However, this algorithm finds an optimal b-coloring of $G$ only when $G$ does not have a good set (i.e., when $\chi_b(G)=m(G)-1$). If $G$ has a good set, we can find an optimal b-coloring of $G$ if a list of all special cycles of $G$ is given. %For instance, if $G$ is an outerplanar graph, then such a list can be found in polynomial time \cite{}.

In \cite{Silva.Maffray.11}, Maffray and Silva conjecture that any graph $G$ with no $K_{2,3}$ as subgraph has b-chromatic number at least $m(G)-1$. Observe that these graphs contain all graphs with girth at least 9; thus, we have given a partial answer to their conjecture. Actually, if their conjecture holds, then $\chi_b\ge m(G)-1$ holds for every $G$ with girth at least 5. However, a different approach is needed as our proof strongly relies on the fact that $\mbox{girth}(G)\geq 9$. Moreover, Theorem \ref{thm:goodset} does not hold for an infinite family of cacti with girth 5, as can be seen in \cite{Campos.etal.09}. This means that the situation where $G$ has no good set is not the only situation where a graph $G$ with girth at least 5 cannot be b-colored with $m(G)$ colors.

\section{Acknowledgements}
This work has been partially supported by Funcap (Fun\-da\c{c}\~ao Cearense de Apoio ao Desenvolvimento Cient\'ifico e Tecnol\'ogico) and CNPq (Conselho Nacional de Desenvolvimento Cient\'ifico e Tecnol\'ogico).

\end{document}